\documentclass{amsart} 
\usepackage{amsfonts}
\usepackage{xcolor}
\usepackage{hyperref}
\hypersetup{colorlinks=true,linkcolor=blue,urlcolor=blue}
\usepackage{graphicx}

\newcommand{\bq}{\begin{quote}}
\newcommand{\eq}{\end{quote}}
\newcommand{\bi}{\begin{itemize}}
\newcommand{\ei}{\end{itemize}}
\newcommand{\bd}{\begin{description}}
\newcommand{\ed}{\end{description}}
\newcommand{\ben}{\begin{enumerate}}
\newcommand{\een}{\end{enumerate}}
\newcommand{\bbm}{\begin{bmatrix}}
\newcommand{\ebm}{\end{bmatrix}}
\newcommand{\bea}{\begin{eqnarray*}}
\newcommand{\eea}{\end{eqnarray*}}

\begin{document}

\title[Unreasonable ineffectiveness of mathematics in biology]{A mathematician's view of the unreasonable ineffectiveness of mathematics in biology.}
\date{06 February 2021}

\author{Alexandre Borovik}
\address{Department of Mathematics, University of Manchester, UK}
\email{alexandre@borovik.net}
\thanks{This is the last pre-publication version of the paper.}

\maketitle

\begin{abstract}
This paper discusses, from a mathematician's point of view, the thesis formulated by Israel Gelfand, one of the greatest mathematicians of the 20th century, and one of the pioneers of mathematical biology:\\
\begin{quote}\small
\emph{ There is only one thing which is more unreasonable than the unreasonable effectiveness of mathematics in physics, and this is the unreasonable ineffectiveness of mathematics in biology.}\\
 \end{quote}

\noindent
\textsc{Disclaimer.} The author writes in his personal capacity and views expressed do not represent position of any other person, corporation, organisation, or institution.
\end{abstract}

\section{Israel Gelfand and his views on the role of mathematics in biology}

Perhaps a disclaimer is necessary: I am a mathematician, not a biologist. I was invited to write this paper because I  found myself in a strange role of a custodian of a particular saying by Israel Gelfand, one of the greatest mathematicians of the 20th century, and a pioneer of mathematical biology. My blog  \cite{borovik-blog} became its principal source:

\begin{quote}\small
Eugene Wigner wrote a famous essay \cite{wigner} on the unreasonable effectiveness of mathematics in natural sciences. He meant physics, of course. There is only one thing which is more unreasonable than the unreasonable effectiveness of mathematics in physics, and this is the unreasonable ineffectiveness of mathematics in biology.
\end{quote}

I wish to confirm that, indeed, I heard these words from Israel Gelfand in  private conversations (and more than once) in about 1995--2005. Beyond that, everything in this paper is my opinion or my reconstruction of Gelfand's view of science and life from my conversations with him; I understand of course that my assessments could be very lopsided.

However, when writing this paper, I located and read papers of a few of Gelfand's earliest collaborators in biology and medicine \cite{Arshavsky,Vasiliev,Vorobiev} and was pleased to discover that my reconstructions were concordant with their memories of him. This gives me  hope that my story contains a reasonable approximation to the truth.

I welcome two papers in this volume, by Blanchard and Longo \cite{Longo} and Rodin \cite{Rodin} which touch on the role of mathematics and biology from perspectives close, but not identical to mine.

I found some further justification of my position in the book \emph{Contemporary Debates in Philosophy of Biology} \cite{Debates} which lists 10 questions and, for each question, contains two papers with completely opposite answers. This paper is an attempt to answer the question
\begin{quote}
\emph{Should we accept Israel Gelfand's assessment of the role of mathematics in biology?}
\end{quote}
And my answer is
\begin{quote}
Yes, we should, for the time being: mathematics is still too weak for playing in biology the role it ought to play.
\end{quote}
I will be happy to see a detailed refutation of my thesis which addresses a number of my concerns raised in the present paper.

Also I think that my stories told here are of general human interest and may be even useful for historians and philosophers of science.

It was not my aim to write any kind of a systematic survey. References are sparse and random and used only as illustrations.

\section{The story starts}

I met Gelfand in 1991  at Rutgers University in the USA, and he immediately dragged me into a research collaboration which lasted for more than decade and was partially summarised in our monograph \cite{CoxeterMatroids}.

Because of Gelfand's peculiar style of work\footnote{See  \cite{topaz}, a brief sketch of Gelfand written by a bemused American mathematician.}, I, although a pure mathematician myself,  was often present during his long conversations with other mathematicians, with mathematical physicists, and with his biologist collaborators, first of all, with Alexander Kister.
  Gelfand's conversations with biologists were mostly about the spacial structure of proteins\footnote{\cite{gelfand1996} is one of the papers produced by Gelfand and Kister in that period. I understand nowadays this type of analysis is heavily computer-based and classified as computational biology.}.

 In  our first conversation Gelfand asked me about my early childhood  mathematical experience, and, specifically, about what   moved me to study mathematics. In my answer I mentioned mathematics correspondence schools which sent to me cute little books on mathematics for schoolchildren, including some books for children written by him.  Gelfand looked at me with suspicion and asked me what I had learned from  his little books. My answer:
  \begin{quote}\small ``the general principle: always start solving a problem by looking at the simplest possible example''
   \end{quote}
 delighted him. This was indeed his principle, he was proud of it, and he systematically applied it throughout all his mathematical work -- but perhaps not in biology: I will return to that later, in Section \ref{sec:adequate}.

 I had never heard the words ``mathematical biology'' from Gelfand -- he always used just ``biology''; in a similar vein, he never used the words ``mathematical physics'' or ``theoretical physics'' -- just ``physics''.

However, Gelfand did a lot of highly nontrivial mathematics and was one of the most influential mathematicians  of the 20th century -- in his thinking, the simplest possible example almost instantly led to very deep mathematics. He also was a mathematical physicist -- and of a very applied kind: for example, he was a mathematical advisor to Andrei Sakharov  in the Soviet H-bomb project and was the head of  the team which carried out critically important calculations \cite[p. 185]{sakharov}); not surprisingly, he had deep knowledge of quantum physics. Gelfand also was one of the pioneers of mathematical biology and  had experience of  50 years of research in that absolutely new,  at his time, area. Sakharov suggests  in his memoirs \cite[p. 219]{sakharov} that the long years of Gelfand's work in mathematical biology may have been motivated by the tragic early death of his son of leukemia (biologists who worked with Gelfand\cite{Vasiliev, Vorobiev} give more detail of this deeply human story).

\section{The controversy and its potential resolution}

I hope I have explained  why Gelfand's remark was  not made off the cuff and deserves some attention. But his view  was contrasted by the \textsc{Wikipedia}\footnote{\textsc{Wikipedia}, \href{https://en.wikipedia.org/wiki/Unreasonable_ineffectiveness_of_mathematics}{Unreasonable ineffectiveness of mathematics}, downloaded 07 Feb 2021.} with the equally strongly expressed opinion of the legendary Leonard Adleman (the `A' in  RSA),  a mathematician, computer scientist, and  cryptographer:

\begin{quote}\small
[In the 1990's] biology was no longer the science of things that smelled funny in refrigerators (my view from undergraduate days in the 1960s [\dots ]). The field was undergoing a revolution and was rapidly acquiring the depth and power previously associated exclusively with the physical sciences. Biology was now the study of information stored in DNA -- strings of four letters: A, T, G, and C and the transformations that information undergoes in the cell. There was mathematics here! \cite[p. 14]{adleman}
\end{quote}

I agree, there is mathematics there. DNA computing, pioneered by Adleman, is a part of mathematics and is fantastic new computer science. But his story is more about application of biology to computer science than application of mathematics to biology. The same could be perhaps be said  about some other recent development, say, the study of  ``artificial life'' \cite{kovitz}.

Also, we have to take into account the fantastic progress of biology, and especially genomics, over the last 20 years  which perhaps makes Gelfand's thesis outdated. It suffices to mention the very recent example: a detailed epigenomic map of non-protein coding segments of human DNA associated with human deseases\footnote{From \href{https://news.mit.edu/2021/epigenomic-map-reveals-circuitry-human-disease-regions-0203}{Epigenomic map reveals circuitry of 30,000 human disease regions, MIT News of February 3, 2021}:
\begin{quote}\footnotesize What we’re delivering is really the circuitry of the human genome. Twenty years later, we not only have the genes, we not only have the noncoding annotations, but we have the modules, the upstream regulators, the downstream targets, the disease variants, and the interpretation of these disease variants
\end{quote}
-- says Manolis Kellis, a professor of computer science, a member of MIT's Computer Science and Artificial Intelligence Laboratory and of the Broad Institute of MIT and Harvard, and the senior author of the new study.} \cite{Boix}.

However,
\begin{itemize}
\item  Gelfand's thesis deserves a discussion. It  should, and can be, discussed  without undermining in any way the successes and heroic efforts of mathematical biologists (Gelfand, after all, was one of them) and bioinformaticians.
\item In his paper, Wigner had in mind pretty highbrow mathematics -- he himself is famous for classifying elementary particles in terms of unitary representations of Lie groups. There is one more  thing which is more unreasonable than the unreasonable effectiveness of ``higher'' mathematics in physics -- and this is the unreasonable effectiveness of arithmetic (even mental arithmetic)  in physics.
 \item The fantastic, explosive  growth of genomics, and studies of RNA and DNA is the evidence for existence of a natural affinity of these parts of biology and mathematics / computer science.
\item But there is more than affinity between mathematics and physics: by their origin, they are twin sisters.
 \item    Biology is much more complex than physics.
\item At its present form, mainstream mathematics  approaches the limits of its
potential applicability to biology. To be useful in the future, mathematics  needs to change dramatically -- and there are good intrinsic reasons for that within mathematics itself.
\end{itemize}

In this paper,  I will try to touch, briefly, on all these points -- but not always in
the same order.

\section{The unreasonable effectiveness of  mental arithmetic  in physics}

It is likely that for Gelfand one of the benchmarks of mathematics' success in applied physics was the creation of the hydrogen bomb -- and he supplied the exceptionally difficult computational part for it. He produced sufficiently precise numerical estimates for processes within the physical event which, most likely, had never  before happened on the surface of the Earth -- radiation implosion. Calculations required digital electronic computers, the first ever -- they were designed and built specifically for that purpose.

But the road to the dawn of the computer era went through tens of thousands of manual (frequently back-of-envelope) calculations and quick mental estimates, say, of physical magnitudes which had never been measured before -- with the aim to get some idea of the size of a measurement device needed and the precision of the measurement required. In physicists' folklore, questions of that kind were known as \emph{Fermi  problems} and could be asked about anything in the world, as Enrico Fermi did, when recruiting young physicists in the Manhattan Project while being unable, for reasons of secrecy, give them any indication of what their future work was about. Instead, he was asking them something like

\begin{quote}\small
How many piano tuners are  in Chicago?
\end{quote}
and invited the interviewees to think aloud, and accessed their reasoning.

Enrico Fermi's  report  \emph{My Observations During the Explosion at Trinity on July 16, 1945}  remains the mother of all mental estimates in physics:

\begin{quote}\small
About $40$ seconds after the explosion, the air blast reached me. I tried to estimate its strength by dropping from about six feet small pieces of paper before, during, and after the passage of the blast wave. Since, at the time, there was no wind I could observe very distinctly and actually measure the displacement of the pieces of paper that were in the process of falling while the blast was passing. The shift was about $2\frac{1}{2}$ meters, which, at the time, I estimated to correspond to the blast that would be produced by ten thousand tons of T.N.T.
\end{quote}

The energy output of the first ever nuclear explosion was calculated, on the spot, and by mental arithmetic, from observation of pieces of paper falling on the ground -- and estimated correctly, as proper measurements confirmed.\footnote{Physicists I spoke to told me they believed that Fermi's calculation was most likely based on the so-called \emph{dimensional analysis} rooted in the scale invariance frequently present in physical phenomena. Kolmogorov's deduction of his ``5/3'' Law (Section \ref{section:Komogorov} and the Appendix)) was also done that way.}

Gelfand was definitely familiar with physicists' love for this kind of mental trick. He told me that he once met Sakharov, who told him: ``You know, on the way here, I did some mental calculation and was surprised to discover that the Sun produces, per unit of mass, less energy  than produced in a pile of rotting manure''\footnote{This observation deserves to be wider known. Life on Earth exists thanks to steady supply of energy from a natural  thermonuclear fusion reactor, safe, clean, stable, reliable, cheap -- our Sun. It is tempting to assume that the promised thermonuclear reactors (already decades in development) will offer the same benefits. But the Sun's power to mass ratio is a bit disappointing. And here is a  Fermi problem for the reader: estimate the size of a pile of manure which would provide an adequate power supply to your home (lights, heating / air conditioning, hot water, all appliances, etc., and add a couple of all-electric cars to the equation), and estimate at what rate the heap has to replenished.}. On hearing this from Gelfand, I  was also surprised and did my own calculations -- Sakharov was (of course)  right. Later I told the story to my astrophysicists friends -- they were astonished, made their calculations (much faster than I did,  I have to admit)  -- and were completely perplexed.

So, this is the way physicists (well, at least experimental physicists) are thinking -- how could it happen to be  so effective?  My proposed answer is in the next section.

\section{Twin sisters: Physics and Mathematics}
\label{twin-sisters}

I will be using the definition (or description) of mathematics given by Davis and Hersh \cite[p. 399]{Davis-Hersh}:

\begin{quote} \small
mathematics is the study of mental objects with reproducible properties.
\end{quote}
The famous mathematician David Mumford uses this formulation
in his paper \cite[p. 199]{Mumford} and further comments on it:

\begin{quote}\small
I love this definition because it doesn't try to limit mathematics to what has been called mathematics in the past
but really attempts to say why certain communications are
classified as math, others as science, others as art, others as
gossip. Thus reproducible properties of the physical world
are science whereas reproducible mental objects are math.
\end{quote}

Mumford's observation can be directly incorporated in (my own) definition:

\begin{quote} \small
mathematics is the study of mental objects and constructions with reproducible properties which imitate the causality structures of the physical world, and are expressed in the human language of social interactions.
\end{quote}

The most basic elements of the causality structures of the world are schemes for expression of observations of the world so self-evident that they never mentioned in physics. For example, if you have some spoons and some forks in your cupboard and you can arrange them in pairs, with no spoon and no fork being singled out, and if you then mix spoons and forks in a box  and start matching them in pairs again, it \emph{must} be a perfect match.

Please notice the word \emph{must} -- its basic use is for expressing relations between people; please also notice that words like `must', `forces', `follows',  `defines', `holds'  etc.\ normally used for description of actions of people and relations between people, play an essential role in any mathematical narrative.\footnote{Without this emphasis on the social interactions language it would be impossible to explain a fact  frequently ignored in discussions of mathematics:  the surprisingly loose and informal language used by mathematicians when they talk about mathematics between themselves -- it has almost nothing in common with the language of published mathematical texts.} What we see in the example with spoons and forks is the mathematical concept of the  one-to-one correspondence between finite sets -- as it appears ``in the wild''.  A mental construction on the top of one-to-one correspondence produces  natural numbers, arithmetic operations, and the order relation. They are interesting for their universal applicability:

\begin{itemize}
\item the number of my children is smaller than the number of protons in the nucleus of Lithium,
\item which, in its turn, is smaller than the number of Galilean moons\footnote{Galilean moons can be objectively defined as satellites of Jupiter visible from Earth via a primitive telescope or standard binoculars.} of Jupiter;
\item which is the same as the  number of  bases of DNA.
\end{itemize}

This is a true statement about four groups of objects in the real world which have absolutely no ``real world'' connections between them.
\begin{quote}
\textbf{The humble natural numbers are already  a huge abstraction}.\footnote{At least one human culture was documented as having no concept of number and no number words in the language: that is of Pirah\~{a} people in the Amazon rainforest \cite[p.260]{Everett}.}
\end{quote}
 The question about the ``unreasonable effectiveness'' has to be asked already about arithmetic, with an obvious answer: yes, arithmetic is effective in biology -- every time we have to count some distinctive and stable objects.

It is a summary of experience accumulated by  humanity over millennia: the causality structures of the physics universe are so robust that their consequences could be developed within mathematics independently from physics -- and remain consistent (that is, do not generate contradictions).

 Moreover, these mathematical developments could happen to be useful for description and modeling of physical phenomena. Ptolemean astronomy was built on the basis of highly developed by that time spherical geometry (born from the needs of astronomy, by the way) in absence of   some key inputs from astronomic observations and from physics which became available only much later -- still, it provided a reasonable approximation to the observed movement of planets in the sky.\footnote{David Khudaverdian kindly explained to me that he does not see any problems with transferring, from the plane to the sphere, of his algorithm (and his computer programme) for approximate reproduction, by a linkage mechanism, and with preservation of the velocity of the point, of movement of a point along a plane curve, see \url{https://david.wf/linkage/theory.html}. It would be interesting to see what this algorithm would do with a kind of data that Ptolemeus could use. This is just a remark on how far we moved from the time of Ptolemeus.}

 At their birth, quantum mechanics and general relativity theory already had their mathematical machinery essentially ready and waiting to be used (perhaps with one important exception, as I'll explain it in minute).  What is important, the efficiency of mathematics in description  the explanation of the real world was demonstrated at least two millennia ago at the level of arithmetic, primitive algebra and geometry. This is a well established historic fact.\footnote{Leonard Adleman was already mentioned here. He is co-inventor of RSA, one of the most widely used cryptographic systems, critically important for the world system of financial transactions, among many other uses. The belief in the security of RSA entirely depends on the assumption that factorisation of integers into products of prime numbers is an exceptionally hard problem. This is  a historic observation extracted from two millennia of human experience with arithmetic. There is still no proof.}

This justifies the motto coined by my colleague Robert A. Wilson:
\begin{quote}\small
Mathematics: solving tomorrow's problems yesterday. \cite{Wilson}
\end{quote}
Of course, occasionally mathematics has to solve today's problems. This had happened with the theory of distributions (or generalised functions): they were invented (or made popular) by one of the founders of quantum physics, Paul Dirac (including the famous $\delta$-function), and were quickly and smoothly incorporated into mathematics; Gelfand was one of the principal contributors to the new theory.

Regarding Gelfand's statement about mathematics and biology, I think he felt  that he faced a much more serious challenge:  the existing mathematics was not directly applicable in biology: some new mathematics was needed. I will return to that point in Section~\ref{sec:adequate}.

And now I wish to offer a mental experiment.

Imagine that over the last 11 thousand years (that is,  the period of stable climate following upon the last ice age  which allowed the  human civilisation to develop) the atmospheric conditions on Earth were different: the skies were always covered, even in the  absence of clouds, by a very  light haze, not preventing the development of agriculture, but   obscuring the stars and turning the sun and the moon into amorphous light spots. Would mathematics have had a chance to develop beyond basic arithmetic and geometry sufficient for measuring fields and keeping records of harvest?  I doubt that. Civilisations which developed serious mathematics also had serious astronomy (it was an equivalent of our theoretical  physics). But I claim even more: the movement of stars in the sky was the  paradigm of precision and reproducibility, the two characteristic features of mathematics. Where else could humans learn the concept of absolute precision?

Speaking about mathematics and physics as twin sisters,  it is almost impossible not to mention the most extreme  point of view on relations between the two sciences.  It belongs to the famous mathematician Vladimir Arnold  \cite{arnold-teaching}:

\begin{quote}\small
\emph{Mathematics is part of Physics.\\
Physics is an experimental discipline, one of the natural sciences.\\
Mathematics is the part of Physics where experiments are cheap.}
\end{quote}

Not every mathematician would agree with that. But I think it is important to put this extreme formulation on record, especially in the context of this paper.

\section{My own doubts about the role of mathematics in biology}

\subsection{My mathematical background}
Everything said in the rest of this paper is my own opinion as a mathematician with 45 years of diverse experiences in mathematics. Over the last 25 years I was engaged -- in parallel with some mainstream and hard core pure mathematics which I was always doing -- with the study of various probabilistic and non-deterministic methods for solving problems in algebra. This made me quite receptive to David Mumford's idea \cite{Mumford} that the future of mathematics is stochastic. I mention this because I believe in the stochastic nature of the underlying laws of biology, whether they are expressed mathematically or not.

This is a huge theme, and in this paper, my arguments are only indicated, not expanded in any detail.

\subsection{Biology as a study of algorithms}
Speaking about  biology, and especially molecular biology, it is important to understand that it is not a natural science in the same sense as physics. It does not study the relatively simple laws of the world. It studies objects which do not exist in physics, and cannot be meaningfully reduced to physical phenomena: \emph{algorithms}.

 It has to deal with molecular algorithms (such as, say, the transcription from DNA to RNA and further translation into synthesis of proteins which ensures the correct spatial shape and the correct functioning of the protein molecule -- and this chain of transformations continues all the way down to specific patterns of neuron firing). Of course I agree with Adleman \cite{adleman}  that  this part:
\begin{quote}\small
``\emph{the transformations that information undergoes in the cell}''
 \end{quote}
can be understood mathematically (or by means of computer science).\footnote{The design of the BioNTech/Pfizer vaccine is enthusiastically greeted by cryptanalysis / computer security geeks who immediately started to ask interesting questions, see Bert Hubert \cite{Hubert} -- but we also should not forget the tremendous work of molecular biologists which made the success possible.  Cryo–electron microscopy resulted in the structure analysis of the SARS-CoV2 spike protein in complex with its cognate cell receptor \cite{Wrapp}, which, in its turn, made possible the design of the stabilized spike protein mutant that has been successfully adapted for the vaccine design for both RNA based  BioNtech/Pfizer and Moderna vaccines.}

Adelman's paper was written in 1998 and stayed within the Central Dogma of molecular biology. He occasionally made even more restrictive statements:
\begin{quote}
  \small\emph{ The polymerase enables DNA to reproduce, which in turn allows cells to reproduce and ultimately allows you to reproduce.  For a strict reductionist, the replication of DNA by DNA polymerase is what life is all about.} \cite[p. 54]{adleman}
\end{quote}
Let us stay for a minute under the umbrella of  the Central Dogma.

 \subsection{Irreversibility}
First of all, we need to take into consideration that there are  many stages of  the transformations ``that information undergoes in the cell'', and  each of them has its own mechanisms for re-encoding the information into a different ``language''.  Each transformation could happen to be a one-way function or procedure, with sufficiently clear ways of performing it, but without rules  for reversing the transformation. Why?

Because all these sophisticated and subtle mechanisms were developed in the course of evolution. The clarity and precision of transformation were obvious selection criteria -- otherwise organisms could not leave  viable  descendants, and, most likely, could not even function themselves.

But it appears that the existence of rules  and mechanisms for reversing each particular transformation had never been a survival criterion. But if some property was not a survival criterion, why we should we expect that it dominates the population? If it was inherited form previous stages of evolution, and lost its usefulness, it is likely that it will be suppressed by mechanism controlling gene expression. (Here we start to deviate from the Central Dogma.)

Avoiding terminology from complexity theory and cryptography, one may still say
\begin{quote}
\textbf{The transformations that information undergoes in the cell is a cascade of functions which could  happen to be not effectively reversible.}
\end{quote}
Without giving a precise definition, I wish to remark that in mathematics such transformations (functions) are called \emph{one-way functions}. A canonical (alleged) example of a one-way function  is multiplication of integers: it is very easy to multiply two integer numbers $p$ and $q$; if  $n = p \times q$, finding factors $p$ and $q$ when given only $n$, is believed to be an impossibly difficult problem, especially if $p$ and $q$ are very large prime numbers. The catch is that it is not proven that factorisation is difficult, it is only a historic observation, the total of experiences accumulated by mathematicians over 2,000 years. The security of the famous RSA cryptosystem, the backbone of electronic finance, is a belief, not a fact.

Almost all mathematicians believe that one-way functions exist, but this remains a conjecture, it is not a theorem, it is not proved. Moreover, almost all functions are likely to be one-way -- but there is no proof of that. On that matter, mathematics is still at square zero.

To summarise,
\begin{quote}
\textbf{Mathematics of nowadays has no tools (and perhaps will never have) for reversing transformations of unknown provenance and of that size of inputs.}
\end{quote}

But, inverting everything that can be reversed is one of the paradigms of mainstream mathematics;  even if you are not a mathematician, recall how many hours you spent at school solving all these equations and systems of equations; this was about it:  reversing   mathematical operations and inverting functions.

Moreover,  more could be said:
\begin{quote}
\textbf{Being understood by some species which would come to existence perhaps hundreds of millions of years later had never been a selection criterion for molecular algorithms  at any stage of their evolution.}
\end{quote}
This basic remark suggests that the current successes of biology is a fantastic achievement which could never have  been taken for granted.

The further we are  from the Central Dogma and the more information transfer paths are discovered in the cell and  between cells, the more complicated and difficult for analysis things become. In particular, if something appears to be reversed, this is not a  full inverse  map -- loss of information is inevitable.  For example, suppression of genes is  an obvious  removal  (perhaps temporary) of no longer needed  information.

\subsection{Black boxes}
My co-author \c{S}\"{u}kr\"{u} Yal\c{c}\i nkaya and I have a toy model for one-way information flows; it is called  \emph{black box algebra} and we are  preparing a monograph on it \cite{BY-monograph}.

Black box algebra  studies categories where objects are some finite mathematical (or computer science, which in this context does not matter) structures called black boxes. Elements of black boxes are binary strings, mathematical operations (perhaps partial) inside of a black box are performed and predicates evaluated  by efficient (in some specific meaning) algorithms. Morphisms are maps from one black box to another which preserve operations and values of predicates in the black boxes  and are performed by efficient algorithms.

There is also a more subtle and flexible relation and subtle relation between black boxes:  we say that a black box $Y$ is\emph{ interpreted} in a black box $X$  if there is an effective map $\alpha: Y \longrightarrow X$ such that for every partial operation (here, for the sake of simplicity of explanation-- binary)  ``$\otimes$''  on $Y$  there is an efficient map $\beta: X \times X \longrightarrow X$ such that
\[
\alpha(y_1 \otimes y_2) = \beta(\alpha(y_1),\alpha(y_2)),
\]
with a  similar property holding for predicates.

The crucial  feature of the theory is that we do not   expect that the inverse morphisms can be also computed efficiently -- morphisms could happen to be one-way functions.  Also, we do not know what is inside of a black box, we can only sample some its random elements and observe their behaviour and interaction with other elements from the sample. In all that there are some analogies with what we see in a cell at a molecular level when we try to look at it from a mathematical point of view.

Black box algebra has  happened to be critically important for solving, by probabilistic methods,  some difficult problems in computational algebra \cite{BY2018,haystack} and is a natural tool for analysis of the so-called homomorphic encryption \cite{BY-homomorphic}.  Some famous intractable problems of algebraic cryptography -- factorisation of integers, the discreet logarithm problem in finite fields and on elliptic curves -- naturally live in the domain of  black box algebra. This shows that this new field of algebra is immensely difficult. This also  supports the nagging feeling that in the world around us almost every process is not reversible (after all, there is  the universal phenomenon of aging followed by the inevitable death), and, moreover, its mathematical description as a function or algorithm (if found)  has no  efficiently computable inverse.

\section{Some further comments on mathematics and evolution}

Once I did some work on genetic (or evolutionary) algorithms in mathematics \cite{BBB,BB}. As it happens in experimental work, not all observations made found their way to publication, especially because my collaborators and I focused on the convergence, in some special cases, of the evolution of a population of non-deterministic algorithms for solving a particular mathematical problem to a known deterministic  algorithm, that is, to an algorithm constructed by humans. For the purpose of this discussion, the cases where the evolution did not progress as we wished would be much more interesting -- and these were the majority of the cases: the geometry of the search space was too complex, and the  evolution of an algorithm  stuck in a  cul-de-sac of a  local maximum.

This raises a natural question: why did  the  evolution of life on Earth produce, and continue to produce, something that apparently works?

Most of the molecular algorithms of life were shaped at the stage of prokaryotes and their immensely complex co-evolution with viruses \cite{Koonin2017,Koonin}. This took, most likely, hundreds of millions  of years, with billions of generations. This number of generations can be reproduced on modern supercomputers. However -- and this was the principal difference from any form of computation that technology might allow us to do -- this was happening in  huge search spaces. The probability of mutations and chances for survival of one of them  in subsequent generations grow with  the size of the population.\footnote{Only very recently, almost a year after the start of the pandemic, I had finally had a chance to  hear  a politician (Shadow Health Secretary in the UK Parliament) referring to  this basic principle in the debate in the British parliament about the pandemic . Still, this is  a colossal success of popularisation of science.}

Any evolution -- an artificial evolution of some artificial entities, or the natural evolution of life -- is blind. In a small number of cases it finds optimal solutions with respect to certain relatively simple constraints and  survival criteria -- the same way as water flows down the slope. For example, all animals living in water, if they have to  be able  to move faster than their prey or predators, have distinctive streamlined shapes dictated by (physical) laws of  fluid dynamics.\footnote{``About 60\%  of the recognized virus taxa have icosahedral capsids, which is unsurprising because the
icosahedron has the largest volume to surface area ratio, closest to that of a sphere, the most thermodynamically favorable three-dimensional shape, and generates the maximum enclosed volume
for shells comprised of a given size subunit.
[\dots ] The other side of the coin, however, is that similar capsid geometries do not necessarily reflect
homologous relationships between viruses: for example, icosahedral capsids emerged at least 11 times
during virus evolution from unrelated CPs with drastically different folds.'' \cite[pp. 4--5]{Koonin}.} However functioning of a cell  means the simultaneous satisfaction of thousands of constraints and criteria. And experiments  show that in problems with multiple constraints evolution does not find an optimal, or even close to optimal, but just a survivable  solution.\footnote{Gregory Cherlin, who read an early draft of my paper, commented at that point:
\begin{quote}\tiny
It is probably looking for solutions to $NP$-complete problems and even with much space is still going to get trapped.  I understand that even the shape of a foam in theoretical physics is a solution to an $NP$-complete problem and nature does not actually produce that shape, even under the laws of physics.
\end{quote}
I share his concern; if $P\ne NP$ (as almost every mathematician expects), mathematisation of biology is likely to be a long slog. Notice that existence of one-way functions implies $P\ne NP$. $P$ vs. $NP$ is one of the Clay Mathematics Institute Millennium Prize Problems, seven problems judged to be among the most important open questions in mathematics.} In short, the surviving solution could be in one of myriads of local optima, sufficiently good to ensure reasonably high probability of survival. Lucky strikes could be so rare that the huge search space and millions of years of evolution produce just one survivable solution, which, as a result, dominates the living world, and is perceived by us as something special.\footnote{Of course, we award this special status, first and foremost, to \emph{ourselves}. There is an almost universal belief that  humanity is the crown of  God's / Evolution's creation.  Ephesians 2:10 is given in the New Living Version as  ``\emph{For we are God's masterpiece}''. The translators of the  (older) King James Version were a bit more modest:  ``\emph{For we are his workmanship}''. It seems that the self-esteem of \emph{H. Sapiens} as a species improves with time.} But it might happen that there is absolutely no reasonably compact external characterization which allows us to distinguish it from other possible solutions, and that its phylogeny (if we will ever know it) is its only explanation.

With the exception of relatively rare periods of regression, evolution progresses  bottom up, from simple to complex. In modern mathematics the situation is different.  Of course, new theories frequently generalise, and are built upon, older theories. However, in concrete research projects and in proving specific theorems mathematics usually works in the opposite direction: from the  more general and abstract down to filling in concrete details.  This is how mathematicians \emph{write} proofs after they got them.  I co-authored a theorem with a proof of 500 pages -- it was published as book \cite{ABC}. Believe me: this could not be done using the ``bottom up'' approach.

The same ``top down'' approach is used in project management: clear identification of priorities and the target, and then planning back to the present position -- with special attention to identification of time critical paths.  The military in more advanced countries  reached real sophistication and efficiency in ``reverse thinking'', both in operational planning and in logistics. In the UK, the army remains the last resort for saving the government's pathetic attempts to manage its response to the COVID epidemic.

I doubt that  the evolution of life had ever done  critical path analysis.

In short, evolution of life  has nothing in common with human problem solving, nothing in common with design and development of mathematical algorithms by mathematicians or computer scientists.

And computers are of no help. I spent considerable time  solving, by non-deterministic methods, mathematical equivalents of the search for the proverbial needle  in a haystack \cite{haystack}.  The biggest structure of that kind where my co-author \c{S}\"{u}kr\"{u} Yal\c{c}\i nkaya and I managed to compute significant and important  substructures, and say something sensible about them, contained about $10^{960}$ elements. The Observable Universe contains around $10^{80}$ electrons. We were computing in something which was $10^{880}$ times bigger than the  Observable Universe. The total number of prokareotes which ever existed on the Earth is nothing in comparison with that. We were successful because we knew what we were looking for, used the powerful global symmetries of the system which we studied, and were able to restrict our work to just  a handful of carefully chosen elements. Also, individual elements were much simpler than any bacteria or archaea; our elements were about $1,000$ bytes long and \emph{had no structure}: we worked with just labels of, or pointers to, random elements -- but bacterial DNA  contains millions of base pairs\footnote{The information content of the  messenger RNA of  BioNTech/Pfeizer vaccine is just above 1 kilobyte.}, and has structure which has to be taken into account if we (humans) try to analyse the DNA molecule. But evolution does not analyse the structure od DNA -- it just checks whether a mutation is advantageous for survival, neutral, or disadvantageous, and these checks are probabilistic by their nature.

And let me repeat:  evolution is blind. 
Evolution does not know what it is looking for. It works via random mutations or exchange of genetic information (again random). For a human mind, even assisted by computers, to navigate the resulting mess -- is a challenging task.

\section{Once more about the unreasonable effectiveness of mathematics in physics}
\label{section:Komogorov}

In the previous Section, non-reversibility of transformations in information flows in cells was highlighted as the principal difficulty for analysing them  mathematically. So it would be useful to look at one of the most extreme cases of the unreasonable effectiveness of mathematics in physics, Andrei Kolmogorov's\footnote{By the way, Gelfand was a student of Kolmogorov.} analysis  of an  incomprehensibly chaotic (and non-reversible, one-way) phenomenon -- \emph{turbulence} --  and try to find: \emph{where is the catch}?

The deduction of Kolmogorov's
seminal ``$5/3$'' Law for the energy distribution in turbulent
fluid \cite{kolmogorov} is so simple that it can be done in a few
lines using only school level algebra (that kind of  derivation can be found in \cite[Section 8.4]{MuM}; I borrow some details from there).

The turbulent flow of a liquid is a cascade of vortices;
the flow in every vortex is made of smaller vortices, all the way
down the scale to the point where the viscosity of the fluid turns
the kinetic energy of motion into heat. So, assume that we are in a steady state, that is, we have a balanced energy flow.

Kolmogorov asked the question: \emph{what is the share of energy carried by vortices of a particular size}?

He got an answer by an elegant short argument based on the important assumption of  \emph{self-similarity} or  \emph{scaling invariance} which  amounted to saying that

\begin{quote}\small
The way bigger vortices are made from smaller ones is the same
throughout the range of wave numbers, from the biggest vortices
(say, like a cyclone covering the whole continent) to a smaller one
(like a whirl of dust on a street corner). \cite{arnold}
\end{quote}
So, this was the catch! And here is Kolmogorov's  formula:
\[
E(k) \approx C \epsilon^{2/3}k^{-5/3}
\]
where $E(k)$  is the \emph{energy density}, $\epsilon$ is \emph{the energy flow},  and $k$ is the \emph{wave number}, while the constant $C$ is dimensionless and is usually close to $1$  (details are in  \cite[Section 8.4]{MuM} or in the Appendix to this paper). 

The status of Kolmogorov's celebrated result is quite remarkable.
In the words of an expert on turbulence, Alexander Chorin \cite{chorin},

\begin{quote} \small
Nothing illustrates better the way in which turbulence is suspended between ignorance and light than the Kolmogorov theory
of turbulence, which is both the cornerstone of what we know and
a mystery that has not been fathomed.
\end{quote}

\begin{quote}\small
The same spectrum [\dots] appears in the sun, in the oceans, and
in man-made machinery. The ``5/3'' Law is well verified experimentally and, by suggesting that not all scales must be computed anew
in each problem, opens the door to practical modelling.
\end{quote}
\normalsize

Vladimir Arnold \cite{arnold} reminds us that the main premises of
Kolmogorov's argument remain unproven---after more than 60 years!
Even worse, Chorin points to the rather disturbing fact that

\begin{quote}\small  Kolmogorov's
spectrum often appears in problems where his assumptions clearly
fail. [\dots] The ``5/3'' Law can now be derived in many ways, often
under assumptions that are antithetical to
Kolmogorov's.  Turbulence theory
finds itself in the odd situation of having to build on its main
result while still struggling to understand it.
\end{quote}

\begin{figure}[h]
\begin{center}
\includegraphics[width=3.6in]{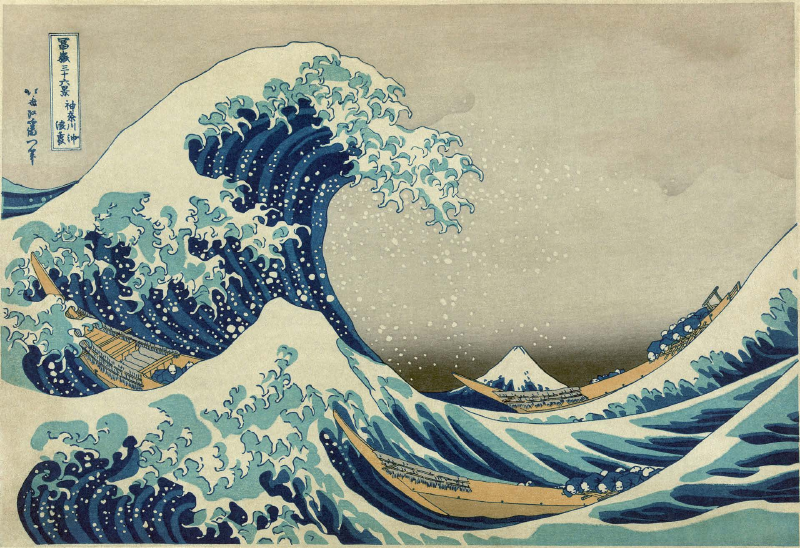}
\end{center}

\caption[\emph{The Great Wave off Kanagawa}  by Katsushika
Hokusai]{{\small Self-similarity and multiple scales in the motion of a fluid, from a
woodcut by Katsushika Hokusai (c.\ 31 October 1760 -- 10 May 1849),  \emph{The Great Wave off Kanagawa}
(from the series \emph{Thirty-six Views of Mount Fuji}, 1823--29).
The distribution of energy between the scales is described by Kolmogorov's ``5/3'' Law.  Luckily for Kolmogorov, water and waves never evolved.  Source:
\emph{Wikimedia Commons}. Public domain.}}

\label{Hokusai}

\end{figure}

This is an interesting case indeed: a remarkable success of mathematics which also shows its limitations. And limitations are obvious: this is only a summarily description of one (although important) aspect of a stochastic phenomenon, Figure~\ref{Hokusai}. In biology, we frequently need something more detailed than that.

A testimony from Sakharov about the role of self-similarity in physics is illuminating:

\begin{quote} \small
Soon after we began the project, I'd come up with an approximate analysis of the important processes specific to the Third Idea\footnote{\emph{The Third Idea} was a more advanced  design of an H Bomb,  much more powerful than the previous versions.}. Mathematically, these were the so-called \emph{self-similar solutions} for certain partial differential equations. [\dots ]
\end{quote}

For Sakharov,  this provided sufficient grounds for gearing up the project:

\begin{quote} \small
 Relying on intuition, and without waiting for
the resolution of
all
theoretical questions or the
final
calculations,
I
issued
instructions and explained to the designers which specifications were critical,
and which could be
adjusted. Through
frequent
visits,
I
established
close
professional
relations with the technical personnel employed
in the design
sector.
I came to appreciate the difficult, painstaking nature of their work, and
the specialized knowledge and talent
it required.
\end{quote}

This is very interesting: an approximate solution based on the assumption of scale invariance  was sufficient for starting the project, but not enough for its completion:

\begin{quote} \small
 Nevertheless, we needed something better than analyses of individual processes using simplified assumptions. Mathematicians at the Installation and in
Moscow worked out new methods for performing complicated calculations by computer. A team headed by Izrail Gelfand, a corresponding member of the Academy of Sciences, played a critical role. I worked closely with Gelfand and his group to develop the basic programs, and we established an excellent rapport despite Gelfand's habit of flying into a rage and shouting
at his colleagues (and sometimes at me as well). After a flare-up, he would stride up
and down his office in silence for a few minutes. When he had regained his
composure, he would return to work and even, on occasion, apologize for his
outburst. Still, I got the impression that Gelfand's colleagues loved him, and
that he had a paternal attitude toward them.
\end{quote}

Basically, Gelfand and his team resolved the extremely difficult  problem numerically, by computer calculations, and without use of the oversimplifying assumption of self-similarity.

This is what I call project management. Importantly, different levels of mathematical modelling were needed at different stages of the development of the project.

At that time, Sakharov was about 35, Gelfand about 45  years old.

\section{Lack of ``global'' scaling invariance in biology}

Molecular level processes within a cell are quite different by their nature from the interaction of cells within a living tissue, and the way an organism is built from its parts is again different. These levels of structural hierarchy developed at different stages of evolution, under different external conditions, and every time evolution had  to use not very suitable means for solving new problems.

When scaling invariance is observed in a living organism, for example, in the form of phyllotaxis \cite{Jean,Lamport,Swinton},  it is  usually restricted to a single level of structural hierarchy.  Not surprisingly, phyllotaxis has happened to be open to mathematical insights, and first serious mathematical study was done by Alan Turing \cite{turing}; he even used one of the first digital electronic computers for related calculations \cite{Swinton}.

The absence of scaling invariance is another obstacle to the effective use of mathematics in biology.

\section{The natural affinity between mathematics and genomics -- and its limits} \label{sec:affinity}

\noindent Returning to the definition of mathematics as ``the study of mental objects with reproducible properties'' (Section \ref{twin-sisters}),
I wish to focus on the word ``reproducible''.\footnote{This section is only a very brief exposition of much more detailed \cite[Chapter 11]{MuM} which contains also analysis of some concrete historic examples.}

\subsection{Memes}
The term \emph{meme} was made popular by Richard Dawkins
\cite{dawkins} and introduced into mainstream philosophy and
cultural studies by Daniel Dennett
\cite{dennett}.  Memes are intended to play
the same role in the explanation of {the} evolution of culture
(and {the} reproduction of individual objects of culture) as genes
{do in the} evolution of life (correspondingly, {the} reproduction
of individual organisms).

The concepts of `meme' and `meme complex' (the latter introduced by Susan Blackmore \cite{Blackmore}) still look
more like  metaphors rather than rigorously defined
scientific terms and have been irreparably undermined by adoption of the word `meme' in social media parlance. In memetics,  specific case studies and
applications (like the one described in \cite{koza}) are still more
interesting than a rather vacuous general theory.

But in discussion of the transmission and reproduction of  mathematics, the meme metaphor has non-trivial aspects.

As I  argue in \cite[Chapter 11]{MuM}, mathematical memes play a crucial role in
many meme complexes of human culture: they increase the precision
of reproduction of the complex, thus giving an evolutionary
advantage to the complex, and, of course, to the memes themselves.
Remarkably, the memes  may remain invisible, unnoticed for
centuries and not recognized as rightly belonging to mathematics.
This  is the characteristic property of
``mathematical'' memes:

\begin{quote}\small If a  meme has the intrinsic property that it increases the
precision of reproduction and error correction of the meme
complexes it belongs to, and if it does that without resorting to
external social or cultural restraints, then it is likely to be an
object or construction of mathematics. \end{quote}

As Ian Stewart put it,

\begin{quote}\small
Mathematics is the ultimate in technology transfer. \cite{stewart}
\end{quote}

Indeed mathematics studies mental objects with reproducible
properties which happen to be built according to highly
precise reproducible rules,
with the precision of reproduction being checked by specific mechanisms, which, in their turn, can also be precisely reproduced and shared. These
rules can themselves be treated as mathematical objects (this is
done in branches of mathematics called mathematical logic and proof theory) and are
governed by metarules, etc. Mathematical objects can reproduce themselves only
because they are built hierarchically. Simple or atomic objects (definitions, formulae, elementary arguments, etc.), form more complicated entities (theorems and their proofs) which, in their turn, are
arranged into theories.

When comparing mathematics with other
cultural systems, we see that some degree of  precision of replication can
usually be found in systems which are relatively simple (like fashion, say). Precision can also be linked to a certain rigidity of the system and an institutionalized resistance to change, as in the case of religion. We do not offer hecatombs to Zeus, but, after 2000 years or so, we still use Euclidean geometry -- and this has happened without anything resembling
the endless religious wars of human history.

Mathematics is so stable as a cultural complex because it has an extremely powerful intrinsic capability for error detection and error correction. The difficulty of explaining the astonishing power of self-correction of mathematics by external factors, social or cultural, is analyzed, in considerable detail, in   \cite{Azzouni}. I claim that the only possible explanation lies in the nature of mathematical memes themselves.

To summarise the role of mathematical objects  in the evolution of human culture, they are memes which happened to be successful and spread because of the following properties:

\begin{itemize}
\item They have extreme resilience and precision of reproduction.

\item  When included in meme complexes (collections of memes which
have better chances for reproduction when present in the genotype as a group), they increase the precision of reproduction of
the complex as a whole.  We will call memes with this property \emph{correctors}.

\item  This error correcting property is intrinsic to mathematics,  its implementation involves only other mathematical objects,  concepts, and procedures --  it does not depend on  external social or cultural restraints.
\end{itemize}

\subsection{Mathematics is huge -- but in comparison with what? } People outside the mathematical community cannot imagine
how big mathematics is. Davis and Hersh point out that between
100 000 and 200 000 new theorems are published every year in mathematical journals around the world. A poem can exist on its own;
although it requires readers who know its language and can understand its allusions, it does not necessarily refer to or quote other poems. A mathematical theorem, as a rule, explicitly refers to other theorems and definitions and, from the instant of its conception
in a mathematician's mind, is integrated into the huge system of
mathematical knowledge. This system remains unified, tightly connected,  and cohesive: if you take two papers at random, say, one on mathematical logic
and one on probability theory, you may easily conclude that they
have nothing in common. However, a closer look at the Mathematics Subject Classification reveals discipline 03B48: Probability and inductive logic.

We see that, despite all this diversity, there is an almost incomprehensible unity of mathematics. It can be compared only with the
diversity and the unity of life. Indeed, all life forms on Earth, in all
their mind-boggling variety, are based on the same mechanisms of
replication of DNA and RNA, and all that genomic stuff looks like mathematics. It is not surprising at all that mathematics and computer science proved to be efficient there. As I have already said earlier, the trouble  with mathematics is likely to start at higher levels of structure of living matter.

Also the comparison with biology is not really in favour of mathematics: it is minuscule in comparison with Life. Allocating, say, 10 kB of  \LaTeX\ code to the proof of a  theorem, 200,000 theorems become 2GB of  \LaTeX\ files. What is 2GB on biology's scale? Nothing. And there is one more issue: texts are only one of media of social transfer of mathematics. A text is alive only while there are people who \emph{wish} and can understand it; alas, their number, per paper, is in single figures. From my experience of a journal editor, I can say that finding a reviewer for a mathematical paper submitted to a journal is becoming increasingly difficult. In the next section I say more about the emerging crisis in mathematics as a cultural and social system.

\subsection{But it looks as if mathematics  is reaching the limits of human comprehension}

Mathematics continues to grow, and if  you look around, you see that mathematical results and concepts involved in practical applications are
much deeper and more abstract and difficult than ever before. And we have to accept
that the mathematics hardwired and coded, say,  in a smartphone, is beyond the reach of the vast majority of graduates from mathematics departments in our universities.

The cutting edge of mathematical research moves further and further away from
the stagnating mathematics education. From the point of view of an aspiring PhD
student, mathematics looks like New York in the \v{C}apek Brothers’ book \emph{A Long Cat Tale} \cite[p. 44]{Capeks}:

\begin{quote}\small
And New York -- well, houses there are so tall that they can’t even finish building them.
Before the bricklayers and tilers climb up them on their ladders, it is noon, so they eat their
lunches and start climbing down again to be in their beds by bedtime. And so it goes on day
after day.
\end{quote}

Joseph and Karel Capek were the people who coined the word `\emph{robot}' for a specific
socio-economic phenomenon: a device or machine whose purpose is to replace a
human worker. Almost a century ago, they were futurologists -- long before the word `futurology' was invented.

Mathematics badly needs its own specialised mathematical robots -- first of all, for checking proofs, which are becoming impossibly long and difficult. One of the more notorious examples -- the Classification of the  Finite Simple Groups (CFSG), one of the central results of the 20th century algebra. In particular, the CSFG underpins quite a number of results and methods in finite combinatorics, critically important for any systematic development of mathematical biology -- after all, no matter how huge they are, protein molecules are built of finitely many atoms.

The original proof of the CFSG, still with holes, was spread over more than 100 journal papers of total length about 15 thousand pages. A proper and structured  proof  is being published, volume by volume, since 1994 \cite{CSFG}. At the present time, 8 volumes out of the originally estimated 12 are published, volume 9 is in print, volume 10 is in preparation, plus 1220 pages of two volumes of an important part of the proof which was developed separately \cite{ASquasithin}. I personally know, I think, almost every person in the world who can read and understand this proof. The youngest of them is Inna Capdeboscq, one of the authors of volume 9; very soon she will be the only non-retired mathematician who understands the proof of the CFSG.

We have to admit that mathematics faces an existential crisis.

Without switching to  systematic use of computer-based proof assistants, and corresponding changes in the way how mathematics is published and taught, mathematics will not be able to face challenges of biology -- moreover, it is likely to enter a spiral of decay.

\section{The search for the adequate mathematical language} \label{sec:adequate}

Israel Gelfand once said to me:

\begin{quote}\small
Many people think that I am slow, almost stupid. Yes, it
takes time for me to understand what people are saying to
me. To understand a mathematical fact, you have to translate it into a mathematical language which you know. Most mathematicians use three, four languages. But I am an old
man and know too many languages. When you tell me something from combinatorics, I have to translate what you say in the languages of representation theory, integral geometry, hypergeometric functions, cohomology, and so on, in too many languages. This takes time. \cite[p. 67]{MuM}
\end{quote}

Gelfand's love to ``simplest possible examples''  as well his insistence on being constantly reminded of the most basic definitions was not a caprice: he used these examples as pointers  toward \emph{ the most adequate mathematical language} for describing  and solving  a particular  problem; if several languages had to be used, he used definitions as synchronisation markers for smooth  translation from one language to another.

I heard from Gelfand these  particular words: ``adequate mathematical language'' many times. I was excited to find the term ``adequate language'' prominently featuring  in reminiscences about him written by  his colleague in neurophysiology, Yuri Arshavsky.

\begin{quote}\small
The widely accepted concept, presently known as the connectionist concept, that the
brain is a form of computational machinery consisting of simple digital elements was
particularly alien to I.M. Gelfand. Everybody in this audience knows that, according to
I.M. Gelfand, the main problem of science is the problem of  ``adequate language.” For a
formulation of adequate logic there must be language that does not simplify a real
situation. His viewpoint was that the situation in which neuroscientists use the language
of electrical spikes and synaptic connections as the only language in their interaction with
the nervous system, should unavoidably lead to principal roadblocks in understanding the
higher, cognitive functions of the brain. Computational models of cognitive functions,
even those looking flawlessly logical and convincing, are usually incorrect, because they
use non-adequate language.  I.M. Gelfand believed that the language of cognitive
neuroscience should be shifted away from the commonly-accepted ``network” language to
the intracellularly-oriented direction. My guess is that this was among reasons for I.M.
Gelfand to shift his biological interests from neurophysiology to cell biology. He used to
ask us –a group of young electrophysiologists, whether we really believed that neurons
do not have, metaphorically speaking, a ``soul,” but only electrical potentials. In other
words, Gelfand’s idea was that the highest levels of the brain include complex, ``smart”
neurons, performing their own functions and that the whole cognitive function is the
result of cooperative work of these complex neurons. As far as I know, most of Gelfand’s
colleagues have been admired by his fantastic intuition in mathematics. I think that
Gelfand’s idea that neurons can have not only electrical potentials, but also a ``soul”
shows that his intuition extended far beyond mathematics. \cite{Arshavsky}
\end{quote}

I strongly recommend  a short Arshavsky's paper \cite{Arshavsky}; in effect it explains

\begin{quote}\small
the unreasonable ineffectiveness of mathematics in neurophysiology
\end{quote}
and explains the shift of Gelfand's interests to cellular biology. The  ``adequate language'' philosophy was not reductionists in the sense that he refused to work within a single structural  layer of living matter. This was his philosophy in mathematics, too.  For example, he insisted that every decent mathematical theory should have a proper combinatorial underpinning (this is why he dragged me into writing a book \cite{CoxeterMatroids} about some exceptionally simple, at first glance, combinatorial objects -- he needed them for his more serious projects).

And, as a fleeting remark: Gelfand's work in medicine also was a quest for an adequate language \cite{Kulikowski}.

Gelfand applied the same ideology to biology. Here, he did not already have a suitable mathematical language at hand -- it had to be developed; perhaps more than one language was needed. The underlying combinatorial theory also did not exist. There was an additional difficulty: unfortunately, in biology simplest possible  examples, which would be natural stating points for this development,  are not so simple. The foundational combinatorics underpinning a description of molecular processes in the cell (of course, if it exists) has to multidimensional -- just look at  the number of degrees of freedom of a large molecule.

I have a feeling that an appropriate multidimensional combinatorics emerges in works of  Alina  Vdovina, see, for example, \cite{Vdovina}, one of several her  works in which she uses, in various contexts, ubiquitous combinatorial structures  made of a  group acting on (or associated in other subtle ways with) a $CW$-complex; at a naive level, these new multidimensional combinatorial structures have rich (local) symmetries and rich and complex branching.

\section{What kind of new mathematics may help?}

Since my paper is not particularly intended  for mathematicians, this section is very brief.

To summarise, my conclusion is that mathematics, as we know it, is unlikely to be effective in biology. We will need to develop some new mathematics for that.

First of all, we need

\begin{itemize}
\item stronger emphasis on stochasticity -- Mumford wrote about that 20 years ago in his paper  \emph{The dawning of the age of stochasticity} \cite{Mumford}, and
\item new multidimensional combinatorics.
\end{itemize}

Also, we need dramatic, fundamental changes in the  everyday work of mathematicians and in the functioning of mathematics as a cultural system; using a biological simile, these have to be changes at the cellular level.

In my opinion, mathematics for biology will be born from the synthesis of three colossal tasks:
\begin{description}
  \item[Stream A] Rebuild  mainstream pure mathematics as a computer based discipline, with routine use of proof assistants and proof checkers (specialist software packages which implement methods of proof theory). Make sure that the use of proof assistants covers all kinds of stochastic stuff and  non-deterministic methods in mathematics.
  \item[Stream B] Introduce methods of AI (artificial intelligence) into computer-based pure mathematics.
  \item[Stream C] Move beyond statistics-based AI, machine learning, data science etc.\footnote{These directions have their share of issues, see, for example, \cite{D'Amour}.} and develop a new kind of AI which  also uses methods of proof theory to provide not only answers, but also structured human-readable explanations and justifications. If necessary, this new AI should be able to generate language and symbolism for these explanations.
\end{description}

The most prominent programme for realisation of Stream A is Vladimir Voevodsky's \emph{Univalent Foundations of Mathematics} \cite{UFM,voevodsky-origins} -- see Andrei Rodin's paper in this volume \cite{Rodin} for discussion of its possible role in biology.  For a very recent example of other developments, see \cite{abraham}. I doubt that the next generation of mathematicians would be willing to handle proofs 500 page long without computer support. Without proof assistants, further progress of mathematics will simply stop, and any talk of mathematics for biology will become meaningless.

A very recent paper \cite{polu} gives a taste of Stream B.

Stream C appears to be a hot ticket in FinTech, with well-funded start-ups (such as \href{http://www.hylomorph-solutions.com/}{Hylomorph Solutions}) fighting for a killer product.

Of course, realisation of this modest proposal will require a dramatic reform of mathematics education (which is dangerous, judging by  grotesque failures of previous attempts).

\section*{Appendix: Kolmogorov's ``5/3''Law}

I borrow this fragment from my book \emph{Mathematics under the Microscope}  \cite[Section.8.4]{MuM}.

The deduction of Kolmogorov'
seminal ``$5/3$'' law for the energy distribution in the turbulent
fluid \cite{kolmogorov} is so simple that it can be done in a few
lines. It  remains the most striking and beautiful example of
dimensional analysis in mathematics.
I was lucky  to study at a good secondary school where my
physics teacher, Anatoly Mikhailovich Trubachov, derived the ``$5/3$'' law in one of his improvised lectures.

 The turbulent flow of a liquid consists of vortices;
the flow in every vortex is made of smaller vortices, all the way
down the scale to the point when the viscosity of the fluid turns
the kinetic energy of motion into heat (Figure~\ref{Hokusai}). If
there is no influx of energy (like the wind whipping up a storm in
Hokusai's woodcut), the energy of the
motion will eventually dissipate and the water will stand still.
So, assume that we have a balanced energy flow,the storm is
already at full strength and stays that way. The motion of a
liquid is made of waves of different lengths;
Kolmogorov asked the question,
what is the share of energy carried by waves of a particular
length?

Here is a somewhat simplified description of his analysis. We
start by making a list of the quantities involved and their
dimensions.

First, we have the \emph{energy flow}(let me recall,
in our setup it is the same as the dissipation of energy). The
dimension of energy is
\[
\frac{\mbox{mass} \cdot \mbox{length}^2}{\mbox{time}^2}
\]
(remember the formula  $K = mv^2/2$ for the kinetic energy of a
moving material point). It will be convenient to make all
calculations \emph{per unit of mass}. Then the energy flow
$\epsilon$ has dimension
\[
\frac{\mbox{energy}}{\mbox{mass}\cdot \mbox{time}} =
\frac{\mbox{length}^2}{\mbox{time}^3}
\]
For counting waves, it is convenient to use the \emph{wave
number}, that is, the number of waves fitting into the unit of
length. Therefore the wave number $k$ has dimension
\[
\frac{1}{\mbox{length}}.
\]
Finally, the \emph{energy spectrum} $E(k)$ is the quantity such
that, given the interval $$\Delta k= k_1-k_2$$ between the two wave
numbers, the energy (per unit of mass) carried by waves in this
interval should be approximately equal to $E(k_1)\Delta k$. Hence
the dimension of $E$ is
\[
\frac{\mbox{energy}}{\mbox{mass}\cdot \mbox{wave number}} =
\frac{\mbox{length}^3}{\mbox{time}^2}.
\]

To make the next crucial calculations, Kolmogorov made the major assumption that amounted to
saying that\footnote{This formulation is a bit cruder than most
experts would accept; I borrow it from Arnold
\cite{arnold}}.

\small\bq The way
bigger vortices are made from smaller ones is the same throughout
the range of wave numbers, from the biggest vortices (say, like a
cyclone covering the whole continent) to a smaller one (like a
whirl of dust on a street corner). \eq\normalsize

Then we can assume that the energy spectrum $E$, the energy flow
$\epsilon$ and the wave number $k$ are linked by an equation which
does not involve anything else. Since the three quantities
involved have completely different dimensions, we can combine them
only by means of an equation of the form
\[
E(k) \approx C \epsilon^x \cdot k^y.
\]
And now the all-important scaling considerations come into the play. In the equation above, $C$ is a constant. Since the  equation should remain the same
for small scale and for global scale events, the shape of the
equation should not depend on the choice of units of measurements,
hence the constant $C$ should be dimensionless.

Let us now check how  the equation looks in  terms of dimensions:
\[
\frac{\mbox{length}^3}{\mbox{time}^2} =
\left(\frac{\mbox{length}^2}{\mbox{time}^3} \right)^x \cdot
\left(\frac{1}{\mbox{length}} \right)^y.
\]
After equating lengths with lengths and times with times and solving the resulting
system of two simultaneous linear equations in $x$ and $y$, we get
\[
x = \frac{2}{3} \;\; \mbox{ and} \;\;  y = -\frac{5}{3}.
\]
Therefore we come to \emph{Kolmogorov's\ ``$5/3$'' Law}:
\[
E(k) \approx C \epsilon^{2/3}k^{-5/3}.
\]
As simple as that.

Basically, I reproduced here the stuff which I first learnt in one of improvised lectures of my  physics teacher at a secondary school, Anatoly Mikhailovich Trubachov -- he derived  the ``$5/3$'' Law as one of examples of usefulness of dimensional analysis.

It is claimed that people like Enrico Fermi, Stanislaw Ulam (co-inventor, with Edward Teller, of the American H Bomb), or Andrei Sakharov, could  do  dimensional analysis off the top of their heads and use it for producing quick on the hoof estimates of various physical quantities or qualitative description of physical processes. By my time  it became a part of mainstream culture in physics -- it could be explained to schoolchildren. 

It is so much simpler than biology\dots

\section*{Acknowledgements}

I thank Andrei Rodin who has encouraged me to write down these my thoughts. I am grateful to  Inna Capdeboscq,  Gr\'{e}goire Cherlin,  Gregory Cherlin,   Nadia Chuzhanova, David Khudaverdian, Zoltan Kocsis,  Alina Vdovina, \c{S}\"{u}kr\"{u} Yal\c{c}\i nkaya, and especially   Anna Voronova,  for their invaluable help.

Two anonymous referees made very helpful and constructive suggestions which improved the paper.

I was  lucky that my biology teacher at the FMSh, the specialist boarding school of the Novosibirsk University, was Zoya Stepanovna Kiseleva; her principal job was research in  molecular biology at the Institute of Cytology and Genetics in Akademgorodok, Novosibirsk. Her lectures focused on molecular biology,  genetics and population genetics. This influenced my much later choice of research problems  in \emph{mathematics}, in particular, my interest in all things nondeterministic. My physics teacher at the same school,

\section*{A comment on bibliography}

This paper is not a systematic survey; bibliographic references are relatively random and serve only for illustrative purposes.


\begin{thebibliography}{99}

\bibitem{abraham} E. \'{A}brah\'{a}m,  J. Davenport, M. England et al.  New \emph{opportunities for the formal proof of
computational real geometry}? arXiv:2004.04034v1 [cs.SC], 2020.

\bibitem{adleman} L. M. Adleman, \emph{Computing with DNA},  Scientific American, August 1998, 34--41.

\bibitem{ASquasithin} M. Aschbacher and Stephen D. Smith.
\textbf{The Classification of Quasithin Groups: I. Structure of Strongly Quasithin K-groups} and \textbf{II. Main Theorems: The Classification of Simple QTKE-groups}. Mathematical Surveys and Monographs, vols. 11 and 112. Amer. Mat. Soc, 2004. (477 + 743 pp.)

\bibitem{ABC} T. Alt\i nel, A. V. Borovik, and G. Cherlin, \textbf{Simple Groups of Finite Morley Rank}, Amer. Math.
Soc. Monographs Series, Amer. Math. Soc., Providence, RI, 2008.

\bibitem{arnold-teaching} V. I. Arnold, \href{https://www.uni-muenster.de/Physik.TP/~munsteg/arnold.html}{ \emph{On teaching mathematics}}, Samizdat (in Russian). Available at many
websites, including \url{https://www.uni-muenster.de/Physik.TP/~munsteg/arnold.html}.

\bibitem{arnold} V. I. Arnold,   \textbf{What is Mathematics?} MTsNMO, Moscow, 2004.

\bibitem{Arshavsky} Yu. I. Arshavsky, \href{https://www.math.rutgers.edu/docman-lister/math-main/events/gelfand-memorial/1309-arshavsky-pdf/file}{Gelfand on mathematics and neurophysiology}, Gelfand Memorial Meetin, Rutgers University,  2009.

\bibitem{Debates} F. J. Ayala and R. Arp, eds. \textbf{Contemporary Debates in Philosophy of Biology}. Wiley, 2009.


\bibitem{Azzouni}  J. Azzouni, \emph{How and why Mathematics is unique as a social practice}. In J. P.  van Bendegem and B. van Kerkhove, eds. \textbf{Perspectives On Mathematical Practices: Bringing Together Philosophy of Mathematics, Sociology of Mathematics, and Mathematics Education},  Springer, 2007.

\bibitem{Blackmore}  S. Blackmore, \textbf{The Meme Machine}. Oxford, Oxford University Press, 1999.

\bibitem{Longo} E. Blanchard and G. Longo, From axiomatic systems to the dogmatic gene and beyond. This volume.

\bibitem{Boix} C. A. Boix et al., \emph{Regulatory genomic circuitry of human disease loci by integrative epigenomics}.  Nature 590, 11 February 2021.

\bibitem{BBB}   R. F. Booth, D. Y. Bormotov, and A. V. Borovik,    \emph{Genetic algorithms and equations in free groups and semigroups}, in \textbf{Computational and Experimental Group Theory}, Contemp. Math., 349 (2004) 63--81.

\bibitem{BB}  R. F. Booth and A. V. Borovik,  \emph{Coevolution of algorithms and deterministic solutions of equations in free groups}, in \textbf{Genetic Programming. 7th European Conference, EuroGP 2004. Coimbra, Portugal, April 2004}  (M. Keijzer et al., eds.) Lecture Notes Comp. Sci. vol. 3003, Springer-Verlag, 2004, pp.~11--22.


\bibitem{MuM} A. V. Borovik, \textbf{Mathematics under the Microscope}. Amer. Math. Soc., Providence, R.I., 2009.

\bibitem{borovik-blog} A. V. Borovik, \href{http://www.borovik.net/selecta/uncategorized/unreasonable-ineffectiveness-of-mathematics-in-biology-2/}{ \emph{Unreasonable ineffectiveness of mathematics in biology}}, Selected Passages From Correspondence With Friends, a blog post,  31 December 2020; \href{https://micromath.wordpress.com/2018/04/14/unreasonable-ineffectiveness-of-mathematics-in-biology/}{an earlier version},  14 April 2018.


\bibitem{CoxeterMatroids} A. V. Borovik,  I. M. Gelfand, and N. White,   \textbf{Coxeter Matroids}, Birkh\"{a}user, 2002.

\bibitem{BY2018} A. Borovik and  \c{S}. Yal\c{c}\i nkaya, \emph{Adjoint representations of black box groups ${\rm PSL}_2(\mathbb{F}_q)$}.  J. Algebra 506 (2018) 540--591. \href{https://doi.org/10.1016/j.jalgebra.2018.02.022} {Available online}.

\bibitem{haystack} A. Borovik and \c{S}. Yal\c{c}\i nkaya, \href{https://www.mub.eps.manchester.ac.uk/in-abstract/adjoint-representations-of-black-box-groups/}{\emph{Searching for a needle in a haystack, which, in its turn, is locked in a big black box}}. In Abstract. The University of Manchester, 2018.

 \bibitem{BY-homomorphic} A. Borovik and \c{S}. Yal\c{c}\i nkaya,   \emph{ Black box algebra and homomorphic encryption}. In: A. M. Bigatti et al., eds. \textbf{Mathematical
Software - ICMS 2020. 7th International Conference. Braunschweig, Germany, July 13–16,
2020}. Lect. Notes Comp. Sci. 12097. Springer 2020. pp. 115--124. arXiv:1709.01169
[math.GR].

\bibitem{BY-monograph} A. Borovik and \c{S}. Yal\c{c}\i nkaya, \textbf{Black Box Algebra}. In preparation.

\bibitem{Capeks} J. Capek and K. Capek, \textbf{A Long Cat Tale}, Prague, Albatros,  1996 (reprinted  from the original of 1927).

\bibitem{chorin} A. J. Chorin, \emph{Book Review: Kolmogorov spectra of turbulence I: Wave turbulence, by V. E. Zakharov, V. S. Lvov, and G. Falkovich}. Bull.
Amer. Math. Soc. 29  no. 2 (1993), 304--306.

\bibitem{D'Amour} A. D'Amour et al. \emph{Underspecification presents challenges for credibility in modern machine learning},  arXiv:2011.03395v2 [cs.LG], 2020.

\bibitem{Davis-Hersh}  P. Davis and R. Hersh, \textbf{The Mathematical Experience}. Boston, Birkh\"{a}user, 1980.

\bibitem{dawkins}  R. Dawkins, \textbf{The Selfish Gene}. Oxford, Oxford University Press,
1976.

\bibitem{dennett}  D. C. Dennett, \emph{Memes and the exploitation of imagination}. The Journal of Aesthetics and Art Criticism 48 (1990) 127--135.

\bibitem{Everett} D. Everett, \textbf{Language: The Cultural Tool}. London,Profile Books, 2012.

\bibitem{gelfand1996} I. M. Gelfand, A. E. Kister, and D. Leshchiner, The invariant system of coordinates of antibody molecules: Prediction of the ``standard'' $C\alpha$ framework of $V_L$ and $V_H$ domains, Proc. Natl. Acad. Sci. USA 93 (April 1996) 3675--3678.

\bibitem{CSFG} D. Gorenstein, R. Lyons, and R. Solomon,    \textbf{The Classification of the Finite Simple Groups}, Mathematical Surveys and Monographs, 40, volumes  1--8 and forthcoming . Amer. Math. Soc., 1994 -- and continued.

    \bibitem{Hersh}  R. Hersh, \textbf{What is Mathematics, Really?} London, Vintage,  1998.

\bibitem{Hubert}     B. Hubert, \href{https://berthub.eu/articles/posts/reverse-engineering-source-code-of-the-biontech-pfizer-vaccine/}{Reverse Engineering the source code of the
BioNTech/Pfizer SARS-CoV-2 Vaccine}, 25 Dec 2020.

\bibitem{Jean} R. V. Jean, \textbf{Phyllotaxis. A systemic study in plant morphogenesis}. Cambridge, Cambridge University Press, 2009.

\bibitem{kolmogorov}    A. N. Kolmogorov, \emph{Local structure of turbulence in an incompressible fluid for very large Reynolds numbers}, Doklady Acad Sci. USSR 31 (1941), 301--305.

\bibitem{kovitz} B.  Kovitz, D. Bender, and M. Poffald, \href{https://www.mitpressjournals.org/doi/abs/10.1162/isal_a_00190}{\emph{Acclivation of virtual fitness landscapes}}, Artificial Life Conference Proceedings, 31 (2019), 380--387.

\bibitem{koza}  J. R. Koza, M. A. Keane and M. J. Streeter, \emph{Evolving inventions}, Scientific American, 288 no. 2 (2003) 44--59.

\bibitem{Koonin2017} M. Krupovic and E. V. Koonin,  \emph{Multiple origins of viral capsid proteins from cellular ancestors}. Proc
Natl Acad Sci USA 114 (2017,) E2401--E2410.

\bibitem{Koonin} M. Krupovic, V. Dolja, and E. Koonin, \emph{Origin of viruses: primordial replicators recruiting capsids from hosts}. Nature Reviews Microbiology, 17 no. 7 (2019), 449--458.

\bibitem{Kulikowski} C. Kulikowski, \href{https://www.math.rutgers.edu/docman-lister/math-main/events/gelfand-memorial/1311-kulikowski-pdf/file}{\emph{Israel Moiseevitch Gelfand and the Search for an Adequate Language for Medical Diagnosis}}, Gelfand Memorial Meetin, Rutgers University,  2009.

\bibitem{Lamport} D. T. A. Lamport, L. Tan, M. Held, and M. J. Kieliszewski,  \emph{Phyllotaxis turns over a new leaf -- A new hypothesis}.
 Int. J. Mol. Sci. 21  (2020), 1145--1159. DOI:10.3390/ijms21031145.



\bibitem{Mumford} D. Mumford, \emph{The dawning of the age of stochasticity}, in \textbf{Mathematics: Frontiers and Perspectives} (V. I. Arnold et al., eds.). Amer.Math. Soc., 2000, pp. 197--218.

\bibitem{polu} S. Polu and I. Sutskever, \emph{Generative language modeling for automated theorem proving}, arXiv:2009.03393v1 [cs.LG], 2020.

\bibitem{Rodin} A. Rodin, Voevodsky's unachieved project. This volume.

\bibitem{sakharov} A. Sakharov,  \textbf{Memoirs}. New York, Alfred A. Knopf, 1990.

\bibitem{stewart}  I. Stewart, \textbf{Does God play dice? The mathematics of chaos}.  London, Penguin, 1990.

\bibitem{Swinton} J.  Swinton, \emph{Watching the Daisies Grow: Turing and Fibonacci Phyllotaxis}. In: Teuscher C. (eds.) \textbf{Alan Turing: Life and Legacy of a Great Thinker}, pp 477--498. Springer, 2004.

\bibitem{Takeuchi} N. Takeuchi and P. Hogeweg,   \emph{Evolutionary dynamics of RNA-like replicator systems: A bioinformatic approach to the origin of life}. Physics of Life Reviews  9 (2012,) 219--263.


\bibitem{topaz} B. I. Topaz, \emph{Conversations with A. S. Golubitski}. Selected Passages From Correspondence With Friends 8 no. 1 (2020), 1--2.

\bibitem{turing} A. M. Turing,  \emph{The chemical basis of morphogenesis}. Phil. Trans. of the Royal Soci. London, B 237 (1952), 37--72.

\bibitem{Vasiliev} Yu. M. Vasiliev, \emph{About I. M. Gelfand's seminar}, Ontogenez 39 no. 6 (2008), 459--461. (In Russian.)

\bibitem{Vdovina} A. Vdovina, \emph{Drinfeld-Manin solutions of the Yang-Baxter equation coming from cube complexes}. arXiv:2007.01163 [math QA].

\bibitem{UFM}  V. Voevodsky,   \href{https://www.math.ias.edu/vladimir/sites/math.ias.edu.vladimir/files/expressions_current.pdf}{Univalent Foundations Project}, 2010.

\bibitem{voevodsky-origins} V. Voevodsky,    \href{https://www.ias.edu/ideas/2014/voevodsky-origins}{The Origins and Motivations of Univalent Foundations}, Institute for Advanced Study, Princeton, 2014.

\bibitem{Vorobiev} A. I. Vorobiev,\emph{ The Biological Seminar of Israel Moiseevich Gelfand}, Ontogenez 39 no. 6 (2008), 462--464. (In Russian.)

\bibitem{wigner} E. Wigner, \emph{The Unreasonable Effectiveness of Mathematics in the Natural Sciences}, Comms.  Pure Applied Maths. 13  no. 1 ( 1960), 1--13.

\bibitem{Wilson} R. A. Wilson, \textbf{Hidden Assumptions},  \url{https://robwilson1.wordpress.com/}.

\bibitem{Wrapp} D. Wrapp et al. \emph{Cryo-EM structure of the 2019-nCoV spike in the
prefusion conformation}. Science 367 no. 6483 (2020), 1260--1263. \url{https://robwilson1.wordpress.com/}.
\end{thebibliography}
\end{document}